\documentclass{ifacconf}

\usepackage{amssymb}
\usepackage{amsfonts}
\usepackage{amsmath}
\usepackage{amsbsy}
\usepackage[dvips]{graphicx}
\usepackage{psfrag}

\usepackage{natbib}

\begin{document}

\begin{frontmatter}
\title{Casimir-Based Control Beyond the Dissipation Obstacle}

\author{Johan Koopman and Dimitri Jeltsema}

\address{Delft Institute of Applied Mathematics, Delft University of
Technology, Mekelweg 4, 2628 CD Delft, The Netherlands}

\begin{abstract}
A prevailing trend in the stabilization of port-Hamiltonian systems
is the assumption that the plant and the controller are both
passive. In the standard approach of control by interconnection
based on the generation of Casimir functions, this assumption leads
to the dissipation obstacle, which essentially means that
dissipation is admissible only on the coordinates of the closed-loop
Hamiltonian that do not require shaping and thus severely restricts
the scope of applications. In this contribution, we show that we can
easily go beyond the dissipation obstacle by allowing the controller
to have a negative semi-definite resistive structure, while
guaranteeing stability of both the closed-loop and the controller.
\end{abstract}


\end{frontmatter}

\section{Introduction}

In the past decade, passivity-based control (PBC) has emerged as a
control design method that respects, and successfully exploits, the
physical structure of a system. Using the port-Hamiltonian (pH)
formalism, state-space control design methods are proposed that lead
to controllers and subsequent closed-loop systems that admit a
physical interpretation. In this approach the Hamiltonian (i.e., the
internal stored energy) of the system is the focal point of the
design method, serving as a Lyapunov function for stability and as a
storage function for passivity; see \citep{Schaft:00} for a basic
introduction and \citep{Duindam:09} for a comprehensive summary of
the developments of the pH framework over the past decade.

Starting with the so-called energy shaping (ES) routine
\citep{Ortega:01}, in which the closed-loop energy is shaped using
static state feedback, numerous extension have led to a variety of
control methods that set out to shape the system's energy,
interconnection and dissipation structure. However, only the static
full-state feedback methods are developed to a form that can be
considered generic. Indeed, the interconnection and
dampingassignment  passivity-based  control (IDA-PBC) method is
shown \citep{Ortega:08} to generate all stabilizing static state
feedback controllers for pH-systems.

Existing dynamic output feedback PBC strategies---the so-called
control-by-interconnection (CbI) methods---are centered around the
notion of Casimir functions, which statically relate the states of
the controller to those of the plant \citep{Schaft:00,Ortega:08}.
Although these Casimir-based control design methods are in some ways
very attractive, they are not yet developed into a generic form.
First of all, the input-output structure of both plant and
controller is assumed to be power-conjugate, i.e., the input and
output of both plant and controller is assumed to constitute a
power-port and are thus of equal dimension. This can be seen as a
drawback, since the input-output structure cannot be chosen
arbitrarily. The second, and more severe drawback is that they are
critically hampered by the so-called dissipation obstacle
\citep{Ortega:01}. The dissipation obstacle dictates that no states
can be stabilized that are subject to pervasive dissipation. Several
methods have been developed to circumvent the dissipation obstacle,
see e.g., \citep{Ortega:01,Ortega:03,Jeltsema:04,Ortega:08}.
However, in general, these methods rely on changing the actual
output of the system and thus require a particular input-output
structure, which is not always possible in practice.

These considerations call for a more general output feedback method.
In this paper, we concentrate on the aforementioned second drawback
and show that by removing the passivity constraint on the
controller, i.e., by allowing for an active controller, the
dissipation obstacle can simply be resolved while stability of the
controller is established by using the Casimir relation between the
plant and controller states. The design of output feedback
controllers for non-collocated input-output systems will be treated
elsewhere.

\section{Control by Interconnection and Casimir Functions}\label{sec:CbI}

In this section, we briefly review the control by interconnection
(CbI) method applied to port-Hamiltonian (pH) plant systems of the
form
\begin{align}\label{eq:pHplant}
\begin{aligned}
\dot{x} &= \big[J(x)-R(x)\big]\frac{\partial H}{\partial x}(x) +G(x)u,\\
y &= G^T(x)\frac{\partial H}{\partial x}(x),
\end{aligned}
\end{align}
where $x\in\mathbb{R}^{n}$, $u,y\in\mathbb{R}^{m}$, $J(x)$ an $n
\times n$ matrix satisfying $J(x)=-J^T(x)$, $R(x)$ an $n \times n$
matrix satisfying $R(x)=R^T(x)$, $G(x)$ an $n \times m$ matrix, and
the Hamiltonian $H(x)$ represents the total stored energy. The
matrix $J(x)$ is usually referred to as the interconnection
structure, while $R(x)$ captures the resistive structure and is
assumed to be positive semi-definite.

If the controller is also a pH system of the form
\begin{align}\label{eq:pHcontroller}
\begin{aligned}
\dot{\xi} &= \big[J_c(\xi)-R_c(\xi)\big]\frac{\partial H_c}{\partial
\xi}(\xi) +G_c(\xi)u_c,\\
y_c &= G^T_c(\xi)\frac{\partial H_c}{\partial \xi}(\xi),
\end{aligned}
\end{align}
where $\xi\in\mathbb{R}^{n_c}$, $u_c,y_c\in\mathbb{R}^{m_c}$,
$J_c(\xi)$ an $n_c \times n_c$ matrix satisfying
$J_c(\xi)=-J^T_c(\xi)$, $R_c(\xi)$ an $n_c \times n_c$ matrix
satisfying $R_c(\xi)=R^T_c(\xi)$, $G_c(\xi)$ an $n_c \times m_c$
matrix, and $H_c(\xi)$ represents the controller energy, then the
interconnection of the plant system (\ref{eq:pHplant}) with
(\ref{eq:pHcontroller}) via the standard (power-preserving) feedback
interconnection $u=-y_c$, $u_c=y$ (assuming $m=m_c$), yields the
closed-loop system
\begin{align}\label{eq:pHcl}
\begin{aligned}
\begin{bmatrix}
\dot x\\
\dot \xi
\end{bmatrix}
&=
\begin{bmatrix}
J(x) - R(x)  & -G(x)G_c^T(\xi)\\
G_c(\xi)G^T(x) & J_c(\xi) - R_c(\xi)
\end{bmatrix}
\begin{bmatrix}
\dfrac{\partial H}{\partial x}(x)\\[1em]
\dfrac{\partial H_c}{\partial \xi}(\xi)
\end{bmatrix},\\
\begin{bmatrix}
y\\
y_c
\end{bmatrix}
&=
\begin{bmatrix}
G^T(x) & 0 \\
0 & G_c^T(\xi)
\end{bmatrix}
\begin{bmatrix}
\dfrac{\partial H}{\partial x}(x)\\[1em]
\dfrac{\partial H_c}{\partial \xi}(\xi)
\end{bmatrix},
\end{aligned}
\end{align}
which is again a pH system.

The closed-loop system (\ref{eq:pHcl}) is stabilized if the
closed-loop Hamiltonian $H(x)+H_c(\xi)$ can be shaped such that it
has a minimum at the desired equilibrium point $x^*$. The usual way
to proceed is by restricting the motion of the closed-loop system to
the subspace
\begin{equation}
\Omega = \big\{(x,\xi)\in\mathbb{R}^{n \times n_c} \big| C(x,\xi)=\kappa \big\},
\end{equation}
with $C(x,\xi)=\xi-S(x)$, where $S(x)$ is assumed to be a differentiable function and $\kappa\in\mathbb{R}$ is some constant,
and such that the closed-loop Hamiltonian becomes
$H(x)+H_c(S(x)+\kappa)$. This is accomplished if, along the
trajectories of (\ref{eq:pHcl}), the functions $S(x)$ are such that
\begin{equation*}
\dot C(x,\xi)\big|_{\Omega}=0.
\end{equation*}
The functions $C(x,\xi)$ are called Casimir functions and are
independent of the Hamiltonian. We are thus looking for solutions
$S(x)$ of the partial differential equations (PDE's)
\begin{align}\label{eq:PDEold}
\begin{bmatrix}
-\dfrac{\partial^T S}{\partial x}(x) & I_{n_c}
\end{bmatrix}
\begin{bmatrix}
J(x) - R(x)  & -G(x)G_c^T(\xi)\\
G_c(\xi)G^T(x) & J_c(\xi) - R_c(\xi)
\end{bmatrix}=0,
\end{align}
which, under the assumption that $R(x) \succeq 0$ and $R_c(\xi)
\succeq 0$, are characterized by the following chain of equalities
\citep{Schaft:00}:
\begin{align}
\dfrac{\partial^T S}{\partial x}(x) J(x) \dfrac{\partial S}{\partial x}(x)
&= J_c(\xi),\label{eq:Jc-old}\\
R(x)\dfrac{\partial S}{\partial x}(x) &= 0, \label{eq:diss-obst} \\[0.5em]
R_c(\xi) &= 0,\label{eq:diss-obst-Rc}\\[0.45em]
J(x)\dfrac{\partial S}{\partial x}(x)  &= -G(x)G_c^T (\xi).
\end{align}

Unfortunately, the application of the CbI method is severely stymied
by the condition (\ref{eq:diss-obst}), which, roughly speaking,
dictates that the Casimir functions cannot depend on the coordinates
that are subject to dissipation. This means that dissipation is
admissible only on the coordinates of the closed-loop Hamiltonian
that do not require shaping. For that reason, this condition is
referred to as the \emph{dissipation obstacle} \citep{Ortega:01}.
However, the dissipation obstacle stems from the assumption that
\emph{both} the plant and controller dissipation structures satisfy
$R(x) \succeq 0$ and $R_c(\xi) \succeq 0$. Although these properties
are necessary  to ensure that the plant and the controller are both
passive systems, they are merely sufficient for passivity of the
closed-loop system. In fact, the passivity assumption of the
controller is unduly restrictive, as is illustrated in the following
example.

\section{A Motivating Example}\label{sec:motexa}

\begin{figure}[h]
\begin{center}
\psfrag{u}[][]{$u$}
\psfrag{p}[][]{$\phi$}
\psfrag{q}[][]{$q$}
\psfrag{C}[][]{$C$}
\psfrag{L}[][]{$L$}
\psfrag{R}[][]{$r$}
\psfrag{+}[][]{$+$}
\psfrag{-}[][]{$-$}
\includegraphics[width=0.3\textwidth]{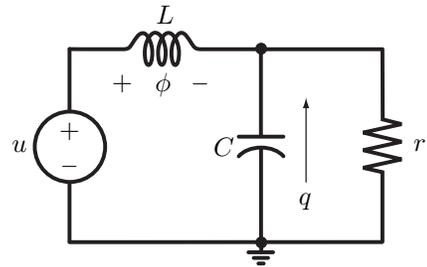}
\caption{RLC with pervasive dissipation.}
\label{fig:RLC}
\end{center}
\end{figure}

Consider the RLC circuit shown in Fig.~\ref{fig:RLC}. Let $\phi$
denote the flux associated to the inductor and $q$ denote the charge
associated to the capacitor, then the equations of motion are given
by the pH description
\begin{align}\label{eq:pH-RLC}
\begin{aligned}
\begin{bmatrix}
\dot \phi\\
\dot q
\end{bmatrix}
&=
\begin{bmatrix}
0 & -1\\
1 & -1/r
\end{bmatrix}
\begin{bmatrix}
\dfrac{\partial H}{\partial \phi}(\phi,q)\\[1em]
\dfrac{\partial H}{\partial q}(\phi,q)
\end{bmatrix}
+
\begin{bmatrix}
1\\
0
\end{bmatrix}u,\\
y&=\begin{bmatrix}
1 & \, 0
\end{bmatrix}
\begin{bmatrix}
\dfrac{\partial H}{\partial \phi}(\phi,q)\\[1em]
\dfrac{\partial H}{\partial q}(\phi,q),
\end{bmatrix},
\end{aligned}
\end{align}
where the Hamiltonian $H(\phi,q)=\frac{1}{2L}\phi^2 +
\frac{1}{2C}q^2$ represents the total stored energy.

This circuit was brought forward in \citep{Ortega:01} as one of the
simplest examples of a system that suffers from the dissipation
obstacle. Indeed, since the equilibrium point equals
$(\phi^*,q^*)=(Lu^*/r,Cu^*)$, which is non-zero for all $u^*
\neq 0$, we need to shape the Hamiltonian in both coordinates, but
condition (\ref{eq:diss-obst}) dictates that
\begin{equation}
\frac{1}{r}\dfrac{\partial S}{\partial q}(\phi,q)=0.
\end{equation}
Hence $S$ can only depend on the flux linkage $\phi$, instead of
both $\phi$ and $q$ as required.

On the other hand, it is easily shown that the circuit is trivially
stabilized by a feedforward control of the form $u=u^*$. Although
this controller is not passive, it does allow for an energy shaping
interpretation since by setting $u=u^*$ the circuit dynamics
(\ref{eq:pH-RLC}) become
\begin{align*}
\begin{aligned}
\begin{bmatrix}
\dot \phi\\
\dot q
\end{bmatrix}
&=
\begin{bmatrix}
0 & -1\\
1 & -1/r
\end{bmatrix}
\begin{bmatrix}
\phi/L\\
q/C
\end{bmatrix}
+
\begin{bmatrix}
1\\
0
\end{bmatrix}u^*\\
&=
\begin{bmatrix}
0 & -1\\
1 & -1/r
\end{bmatrix}
\begin{bmatrix}
\phi/L\\
q/C
\end{bmatrix}
+
\begin{bmatrix}
0 & -1\\
1 & -1/r
\end{bmatrix}
\begin{bmatrix}
0 & -1\\
1 & -1/r
\end{bmatrix}^{-1}
\begin{bmatrix}
1\\
0
\end{bmatrix}u^*\\
&=
\begin{bmatrix}
0 & -1\\
1 & -1/r
\end{bmatrix}
\begin{bmatrix}
\phi/L - u^*/r\\
q/C - u^*
\end{bmatrix},
\end{aligned}
\end{align*}
which can be associated to a shaped Hamiltonian of the form
\begin{equation*}
H_d(\phi,q)=\frac{1}{2L}\left(\phi-Lu^*/r\right)^2 + \frac{1}{2C}\left(q-Cu^*\right)^2
\end{equation*}
such that
\begin{align*}
\begin{bmatrix}
\dot \phi\\
\dot q
\end{bmatrix}
&=
\begin{bmatrix}
0 & -1\\
1 & -1/r
\end{bmatrix}
\begin{bmatrix}
\dfrac{\partial H_d}{\partial \phi}(\phi,q)\\[1em]
\dfrac{\partial H_d}{\partial q}(\phi,q)
\end{bmatrix}.
\end{align*}

\section{Beyond the Dissipation Obstacle}

In the previous section, it is shown that the dissipation obstacle
is a direct consequence of the passivity requirement on both the
system and the controller. In this section, we show that by removing
this constraint, the dissipation obstacle is easily resolved. The
issue that remains then, is how the closed loop can be rendered
stable. Using the Casimir relation between the controller state
$\xi$ and the plant state $x$, we show that the static state
feedback interpretation of the closed-loop system solves this issue
in a straightforward manner.

Consider again the system of PDE's (\ref{eq:PDEold}), or equivalently,
\begin{align}
\dfrac{\partial^T S}{\partial x}(x) \big[J(x)-R(x)\big] - G_c(\xi)G^T(x)
&= 0,\label{eq:PDE1}\\
\dfrac{\partial^T S}{\partial x}(x) G(x)G_c^T(\xi) + \big[J_c(\xi) -
R_c(\xi)\big] &= 0\label{eq:PDE2}.
\end{align}
Substituting (\ref{eq:PDE1}) into the transposed of (\ref{eq:PDE2})
yields
\begin{align*}
 \dfrac{\partial^T S}{\partial x}(x) \big[J(x)-R(x)\big]
 \dfrac{\partial S}{\partial x}(x) - \big[J_c(\xi) + R_c(\xi)\big] &= 0,
\end{align*}
which, after separation of the symmetric and skew-symmetric part,
implies for the relationship between the controller and the plant
dissipation and interconnection structure that
\begin{align}
R_c(\xi) &= -\dfrac{\partial^T S}{\partial x}(x) R(x)
\dfrac{\partial S}{\partial x}(x)\label{eq:Rc}\\
J_c(\xi) &= \dfrac{\partial^T S}{\partial x}(x) J(x) \dfrac{\partial
S}{\partial x}(x)\label{eq:Jc},
\end{align}
respectively. It is directly observed that (\ref{eq:Jc}) coincides
with condition (\ref{eq:Jc-old}), whereas (\ref{eq:Rc}) coincides
with (\ref{eq:diss-obst}) and (\ref{eq:diss-obst-Rc}) if, and only
if, we only allow for dissipation structures  satisfying $R(x)
\succeq 0$ and $R_c(\xi) \succeq 0$. However, if we also allow for a
non-positive controller dissipation structure, we can simply proceed from the
closed-loop plant dynamics
\begin{equation*}
\dot x = \big[J(x)-R(x)\big]\dfrac{\partial H}{\partial x}(x)
- G(x)G_c^T(\xi)\dfrac{\partial H_c}{\partial \xi}(\xi),
\end{equation*}
which, by transposing (\ref{eq:PDE1}), i.e.,
\begin{align*}
G(x)G_c^T(\xi) = - \big[J(x) + R(x)\big]\dfrac{\partial S}{\partial x}(x),
\end{align*}
can be rewritten as
\begin{align}
\dot x = \big[J(x) &- R(x)\big]\dfrac{\partial H}{\partial x}(x)\nonumber\\
& + \big[J(x) + R(x)\big]\dfrac{\partial S}{\partial x}(x)
\dfrac{\partial H_c}{\partial \xi}(\xi).\label{eq:plant-closed-loop}
\end{align}
Hence, if $J(x)-R(x)$ is invertible, we can write
\begin{align}\label{eq:ES_cl}
\dot x =& \big[J(x)-R(x)\big]\left(\dfrac{\partial H}{\partial x}(x)
+ \big[J(x)-R(x)\big]^{-1} \right. \nonumber \\
& \quad  \left. \times \big[J(x) + R(x)\big]
\dfrac{\partial S}{\partial x}(x)
\dfrac{\partial H_c}{\partial S}(S(x)+\kappa)\right),
\end{align}
which, after suitable choices of the controller Hamiltonian
$H_c(S(x)+\kappa)$, may be interpreted as an energy shaping (ES) process
such that the closed-loop plant dynamics take the form
\begin{equation}\label{eq:ES-interpretation}
\dot x = \big[J(x)-R(x)\big]\dfrac{\partial H_d}{\partial x}(x),
\end{equation}
with
\begin{align*}
\dfrac{\partial H_d}{\partial x}(x) & = \dfrac{\partial
H}{\partial x}(x)
+ \big[J(x)-R(x)\big]^{-1}\\
& \times \big[J(x) + R(x)\big]\dfrac{\partial S}{\partial
x}(x) \dfrac{\partial H_c}{\partial S}(S(x)+\kappa)\nonumber,
\end{align*}
satisfying
\begin{equation}\label{eq:Poincare}
\dfrac{\partial^2 H_d}{\partial x^2}(x)= \left(\dfrac{\partial^2 H_d}{\partial x^2}(x)\right)^T\!\!\!\!, \ \text{for all $x \in \mathbb{R}^n$.}
\end{equation}
Furthermore, if $R(x)\succeq 0$, and
\begin{equation}\label{eq:sec-der-test}
\dfrac{\partial H_d}{\partial x}(x^*) \equiv 0, \quad \dfrac{\partial^2 H_d}{\partial x^2}(x^*) \succ 0,
\end{equation}
it follows that $x^*$ is a stable equilibrium of
(\ref{eq:plant-closed-loop}).

In general, stability of (\ref{eq:plant-closed-loop}) does not necessarily
imply closed-loop stability. However, since $S:\mathbb{R}^n \to \mathbb{R}^{n_c}$ is  differentiable by assumption (as is needed in (\ref{eq:PDEold})), it is continuous. By virtue of the well-known preservation of convergence under continuous mappings \citep{Kolmogorov:99}, as $x$ converges to $x^*$, the controller state $\xi$ converges to $\xi^*=S(x^*)$. Hence, stability of $x$ implies stability of $\xi$. This means that the controller
does not need to have a positive semi-definite dissipation structure
in order to be  stable in the closed loop.

Although the method described above extends the traditional
CbI method considerably, it is still hampered by the assumption that
$J(x)-R(x)$ is invertible, and that (\ref{eq:Poincare}) should
be satisfied. Both these assumptions are needed for a pure ES
interpretation of this CbI method, but are, in general, overly
restrictive. However, using arguments similar to those of the static
state feedback IDA-PBC method, the CbI method above can be extended
to a dynamic output feedback IDA-PBC method. Indeed, starting from (\ref{eq:plant-closed-loop}), we then have to look for matrices $J_d(x)=-J_d^T(x)$ and $R_d(x)=R_d^T(x)$, and a function $H_d:\mathbb{R}^n \to \mathbb{R}$ satisfying (\ref{eq:Poincare}), such that 
\begin{align}
\dot{x} &= \big[J(x) - R(x)\big]\dfrac{\partial H}{\partial x}(x)\nonumber\\ 
& \qquad + \big[J(x)+  R(x)\big]\dfrac{\partial S}{\partial x}(x) \dfrac{\partial H_c}{\partial S}(S(x)+\kappa)\nonumber\\
    &\equiv \big[J_d(x)-R_d(x)\big]\dfrac{\partial H_d}{\partial x}(x).\label{eq:IDA-interp}
  \end{align}
This result effectively generalizes both the standard passive ES Casimir-based control method as presented in Section \ref{sec:CbI}, and the extended ES Casimir-based control method of (\ref{eq:plant-closed-loop})--(\ref{eq:ES-interpretation}), to a dynamic output feedback IDA-PBC strategy. Summarizing, we have the following proposition.   

\bigskip

\hrule
{\it Proposition.} Consider the interconnection of the plant (\ref{eq:pHplant}) with the controller (\ref{eq:pHcontroller}), assume that (\ref{eq:PDEold}) holds, and that the closed-loop plant dynamics satisfy (\ref{eq:IDA-interp}) (resp.~(\ref{eq:ES-interpretation})). Then, if the closed-loop plant Hamiltonian $H_d(x)$ satisfies (\ref{eq:sec-der-test}) and $R_d(x)\succeq 0$ (resp.~$R(x) \succeq 0$), the equilibrium point $(x^*,S(x^*))$ of the overall system (\ref{eq:pHcl}) is stable.

\medskip
 
\hrule

\section{A Motivating Example (Cont'd)}

Let us return to the RLC circuit of Fig.~\ref{fig:RLC}. In
Section \ref{sec:motexa}, we have seen that the circuit suffers from
the dissipation obstacle. However, suppose we interconnect, again
via standard feedback, the circuit (\ref{eq:pH-RLC}) with a pH
controller of the form
\begin{align}\label{eq:pH-RLCcontroller}
\begin{aligned}
\dot{\xi} &= -R_c\frac{\partial H_c}{\partial \xi}(\xi) +G_cu_c,\\
y_c &= G_c\frac{\partial H_c}{\partial \xi}(\xi).
\end{aligned}
\end{align}
Then, proceeding from (\ref{eq:PDEold}), we find
\begin{equation*}
\left.\begin{aligned}
\dfrac{\partial S}{\partial q}(\phi,q) - G_c &= 0\\
\dfrac{\partial S}{\partial \phi}(\phi,q) + \frac{1}{r}
\dfrac{\partial S}{\partial q}(\phi,q) &= 0\\
\dfrac{\partial S}{\partial \phi}(\phi,q)G_c - R_c &=0
\end{aligned}\right\} \ \Leftrightarrow \ S(\phi,q) = G_c\left(q-\phi/r\right),
\end{equation*}
and $R_c = -\frac{1}{r}G_c^2$, which is non-positive for all
$r<\infty$ and $G_c \neq 0$.

\subsection{Recovering the Feedforward Controller}

Now, setting the controller Hamiltonian $H_c(\xi) = \xi$, with
$\xi=S(\phi,q)$, the dynamics of the controlled circuit becomes
\begin{align}\label{eq:pH-RLC-controlled1}
\begin{aligned}
\begin{bmatrix}
\dot \phi\\
\dot q
\end{bmatrix}
&=
\begin{bmatrix}
0 & -1\\
1 & -1/r
\end{bmatrix}
\begin{bmatrix}
\dfrac{\partial H}{\partial \phi}(\phi,q)\\[1em]
\dfrac{\partial H}{\partial q}(\phi,q)
\end{bmatrix}
+
\begin{bmatrix}
1\\
0
\end{bmatrix}(-y_c)\\
&=
\begin{bmatrix}
0 & -1\\
1 & -1/r
\end{bmatrix}
\begin{bmatrix}
\dfrac{\partial H}{\partial \phi}(\phi,q)\\[1em]
\dfrac{\partial H}{\partial q}(\phi,q)
\end{bmatrix}\\
& \qquad
-
\begin{bmatrix}
0 & -1\\
1 & -1/r
\end{bmatrix}
\begin{bmatrix}
0 & -1\\
1 & -1/r
\end{bmatrix}^{-1}
\begin{bmatrix}
1\\
0
\end{bmatrix}G_c\frac{\partial H_c}{\partial \xi}(\xi)\\
&=
\begin{bmatrix}
0 & -1\\
1 & -1/r
\end{bmatrix}
\begin{bmatrix}
\dfrac{\partial H}{\partial \phi}(\phi,q) + \dfrac{G_c}{r}
\dfrac{\partial H_c}{\partial \xi}(\xi)\\[1em]
\dfrac{\partial H}{\partial q}(\phi,q) + G_c
\dfrac{\partial H_c}{\partial \xi}(\xi)
\end{bmatrix}.
\end{aligned}
\end{align}
Since $\frac{\partial H_c}{\partial \xi}(\xi)=1$, the stabilization
problem is trivially solved by setting $G_c=-u^*$ such that
(\ref{eq:pH-RLC-controlled1}) is equivalent to
\begin{align}\label{eq:pH-RLC-controlled2}
\begin{aligned}
\begin{bmatrix}
\dot \phi\\
\dot q
\end{bmatrix}
=
\begin{bmatrix}
0 & -1\\
1 & -1/r
\end{bmatrix}
\begin{bmatrix}
\dfrac{\partial H_d}{\partial \phi}(\phi,q)\\[1em]
\dfrac{\partial H_d}{\partial q}(\phi,q)
\end{bmatrix},
\end{aligned}
\end{align}
with shaped Hamiltonian
\begin{equation*}
  H_d(\phi,q)=\frac{1}{2L}\left(\phi-Lu^*/r\right)^2 +
  \frac{1}{2C}(q-Cu^*)^2.
\end{equation*}

At this point it is important to emphasize that, although the
circuit is stabilized by a constant control, the underlying
controller dynamics (\ref{eq:pH-RLCcontroller}) take the form
\begin{align}\label{eq:pH-RLCcontroller2}
\begin{aligned}
\dot{\xi} &= \frac{1}{r}(u^*)^2 - u^* u_c,\\
y_c &= -u^*.
\end{aligned}
\end{align}
Despite the fact that $R_c < 0$, for all $u^* \neq 0$, the
controller dynamics is asymptotically stable since from
(\ref{eq:pH-RLC-controlled2}) we may deduce that $\phi \to
Lu^*/r$ and  $q \to Cu^*$, as $t \to \infty$, which, since
$\xi=S(\phi,q)$, implies that the controller state converges to
$\xi^*=-u^*(q^* - \phi^*/r)$. Also note that $u_c \to u^*/r$ implies
$\dot{\xi} \to 0$.

The Casimir approach thus leads us naturally to the feedforward
controller that is known to stabilize the system, and in spite of
its feedforward character, the controller dynamics are instrumental
for the construction of the controller.

\subsection{Output Feedback Control}

Let us next exploit the freedom in choosing $H_c(\xi)$ as
\begin{align*}
  H_c(\xi) &= \frac{1}{2} a_1 \xi^2 + a_2 \xi,
\end{align*}
with $\xi=S(\phi,q)$ and $a_i\in\mathbb{R}$, for $i=1,2$, some
constants to be defined. With this choice, however, it is not
possible to interpret the control action as an ES
process (\ref{eq:ES-interpretation}) since
\begin{align*}
\begin{bmatrix}
\dfrac{\partial H}{\partial \phi}(\phi,q)\\[1em]
\dfrac{\partial H}{\partial q}(\phi,q)
\end{bmatrix}
+
\begin{bmatrix}
1 & 2/r\\
0 & 1
\end{bmatrix}
\begin{bmatrix}
\dfrac{\partial S}{\partial \phi}(\phi,q)\\[1em]
\dfrac{\partial S}{\partial q}(\phi,q)
\end{bmatrix}\dfrac{\partial H_c}{\partial S}(S(\phi,q))
\end{align*}
is not a gradient vector field to which we can associate a
Hamiltonian $H_d(\phi,q)$ satisfying (\ref{eq:Poincare}). 
On the other hand, the closed-loop dynamics of the plant take the form
\begin{align*}
\begin{bmatrix}
\dot \phi\\
\dot q
\end{bmatrix}
&=
\begin{bmatrix}
a_1G_cL/r & -(1+a_1G_c C)\\
1 & -1/r
\end{bmatrix}
\begin{bmatrix}
\phi/L\\
q/C
\end{bmatrix}+
\begin{bmatrix}
-G_c a_2\\
0
\end{bmatrix},
\end{align*}
which can be rewritten in pH form as
\begin{align}
\begin{aligned}
\begin{bmatrix}
\dot \phi\\
\dot q
\end{bmatrix}
&= \Bigg(
\underbrace{\begin{bmatrix}
0 & -1-a_1G_cC/2\\
1+a_1G_cC/2 & 0
\end{bmatrix}}_{J_d}\\
& \qquad
-
\underbrace{\begin{bmatrix}
-a_1G_cL/r & a_1G_cC/2\\
 a_1G_cC/2 & 1/r
\end{bmatrix}}_{R_d}
\Bigg)
\begin{bmatrix}
\dfrac{\partial H_d}{\partial \phi}(\phi,q)\\[1em]
\dfrac{\partial H_d}{\partial q}(\phi,q)
\end{bmatrix},
\end{aligned}
\end{align}
where
\begin{equation*}
H_d(\phi,q)=\frac{1}{2L}\left(\phi-L\alpha/r\right)^2 + \frac{1}{2C}(q-C\alpha)^2,
\end{equation*}
with $\alpha=-r^2 a_2G_c(r^2+a_1CG_cr^2-a_1G_cL)^{-1}$.

Note that the closed-loop plant dynamics is stable if the controller parameters $a_1$ and $G_c$ are selected such that $R_d \succeq 0$. The desired equilibrium point is determined by an appropriate selection of $a_2$.



\section{Final Remarks and Outlook}

In this paper, an extension of the Casimir-based
control-by-interconnection method is presented that resolves the
so-called dissipation obstacle. This result is based on removing the
passivity constraint on the controller and using the Casimir
relation between the plant and controller states to guarantee
stability of the closed loop.

Using a simple RLC circuit example with pervasive dissipation, two
intriguing controller synthesis solutions are presented. First, it
is shown that this dynamic output-feedback based design methodology
is able to produce a feedforward controller, while the second
solution shows how the dynamics of the controller are instrumental
in shaping the energy in the coordinates that are not used for
feedback.

Although the results presented in this paper only emphasize the
possibility of an active controller, future work will include
conditions under which also active plants are allowed, as well as
the consideration of pH systems with direct feedthrough and
output feedback for non-collocated input-output channels.

\newpage

\bibliographystyle{asmems4}

\end{document}